\documentclass[12pt]{article}

\usepackage{amssymb,amsthm,amsmath}
\usepackage[usenames,dvipsnames]{xcolor}
\newcommand{\hwp}{\mathrm{HWP}}
\usepackage{pgf,tikz}
\allowdisplaybreaks
\usepackage{color,graphicx,epsfig}
\usepackage{mathrsfs}

\newcommand{\comment}[1]{}
\definecolor{teal}{RGB}{0,128,128}
\definecolor{darkpurple}{RGB}{128,0,128}

\usepackage[titletoc,toc,title]{appendix}

\newtheorem{theorem}{Theorem}[section]

\newtheorem{lemma}[theorem]{Lemma}
\newtheorem{cor}[theorem]{Corollary}

\theoremstyle{definition}
\newtheorem{defn}[theorem]{Definition}

\theoremstyle{definition}
\newtheorem{definition}[theorem]{Definition}

\def \cF {{\cal F}}

\def \Z {\mathbb Z}

\title {
A survey on constructive methods for the Oberwolfach problem and its variants
}

\author{
A.\ C.\ Burgess \footnotemark[1]  
\and 
P.\ Danziger    \footnotemark[2] 
\and
T.\ Traetta     \footnotemark[3] 
}

\begin{document}
\date{}

\maketitle

\footnotetext[1]{Department of Mathematics and Statistics, University of New Brunswick, 100 Tucker Park Rd., Saint John, NB  E2L 4L5, Canada. Email: andrea.burgess@unb.ca}
\footnotetext[2]{Department of Mathematics, Toronto Metropolitan University, 350 Victoria St., Toronto, ON  M5B 2K3, Canada. Email: danziger@ryerson.ca}
\footnotetext[3]{DICATAM, Universit\`{a} degli Studi di Brescia, Via Branze 43, 25123 Brescia, Italy. E-mail: tommaso.traetta@unibs.it}

\abstract{
	The generalized Oberwolfach problem asks for a decomposition of a graph $G$ into specified 2-regular spanning subgraphs $F_1,\ldots, F_k$, called factors. 
	The classic Oberwolfach problem corresponds to the case when all of the factors are pairwise isomorphic, and $G$ is the complete graph of odd order or the complete graph of even order with the edges of a $1$-factor removed. 
	When there are two possible factor types, it is called the Hamilton-Waterloo problem.
	In this paper we present a survey of constructive methods which have allowed recent progress in this area. 
	Specifically, we consider blow-up type constructions, particularly as applied to the case when each factor consists of cycles of the same length. 
	We consider the case when the factors are all bipartite (and hence consist of even cycles) and a method for using circulant graphs to find solutions.
	We  also consider constructions which yield solutions with well-behaved automorphisms.
	}

\section{Introduction}

The Oberwolfach problem was famously introduced in 1967 by Gerhard Ringel at a conference in Oberwolfach Germany. 
At this conference it was the tradition for attendees to have dinner together each night, Ringel asked if it was possible to sit the attendees at the round dinner tables over the course of the conference so that each attendee sat next to each other attendee exactly once. 
As such he was asking for a 2-factorization of the complete graph into a specified 2-factor. 

More specifically, given a graph $G$ with vertex set $V(G)$, a factor $F$ is a spanning subgraph of $G$, that is, $V(F)=V(G)$. We refer to a $k$-factor when every vertex of $F$ has degree $k$. So, a 1-factor is a perfect matching and a 2-factor is a collection of disjoint cycles (a $2$-regular graph) which span the vertices of $G$. 
Letting $2\leq \ell_1 <\ell_2< \ldots< \ell_t$, we denote by 
$[\,^{\alpha_1}\ell_1,\, ^{\alpha_2}\ell_2,\, \ldots, \,^{\alpha_t}\ell_t]$ any $2$-regular graph containing $\alpha_i$ cycles of length $\ell_i$, for $i=1,2, \ldots, t$. If all cycles of $F$ have the same length, we speak of a uniform $2$-factor.
A factorization of $G$ is a collection $\mathcal{F}$ of factors whose edge sets partition between them the edge set $E(G)$ of $G$. 
If all factors of $\mathcal{F}$ are $k$-factors (resp. isomorphic to $F$), we refer to $\mathcal{F}$ as a $k$-factorization (resp. an $F$-factorization) of $G$. 

The Oberwolfach problem can be viewed as a generalization of the classical problem of factoring $K_v$ into triples, first proposed by Kirkman \cite{K50} in 1850.
This problem was independently solved 
by Jiaxi Lu \cite{Lu}
and by Ray-Chadhuri and Wilson \cite{RCW}, more than a century later.

Whereas classical design theory views triples as $K_3$ and generalizes to ask for  factorizations of $K_v$ into complete graphs, the Oberwolfach problem views triples as 3-cycles, and generalizes to decomposing $K_v$ into 2-factors. 
More generally, given a graph $G$ of order $v$ and a multiset $\mathcal{F}=\{\,^{\alpha_1}F_1, \ldots, \,^{\alpha_t}F_t\}$ of $2$-factors of $G$, where each $\alpha_i$ is positive and
the factors $F_i$ are pairwise non-isomorphic,
the Generalized Oberwolfach Problem GOP$(G; \cF)$ asks for a $2$-factorization 
$\mathcal{F}' = \left\{F'_{ij} \mid 1\leq i \leq t, 1\leq j \leq \alpha_i\right\}$ of $G$ where each $F'_{ij}$ is isomorphic to $H_i$.  Necessarily, $G$ must be regular of degree $2\sum_i \alpha_i$. When each $F_i$ is uniform with cycles of length $\ell_i$, in the GOP notation we may replace $\mathcal{F}$ with $\{\,^{\alpha_1}\ell_1, \ldots, \,^{\alpha_t}\ell_t\}$

For $t=1$, the problem is referred to as the Oberwolfach problem (on the graph $G$)
and denoted by $OP(G; F)$. We reserve the notation GOP to the case where $t\geq 2$.
We point out that when $t=2$, one usually speaks of the Hamilton-Waterloo problem 
(on the graph $G$) denoted by HWP$(G; \,^{\alpha_1}H_1, \,^{\alpha_2}H_2\}$. 

We mainly focus on GOP$(G; \cF)$ for the following classes of graphs $G$:
\begin{enumerate}
  \item $G=K_{2n+1}$ is the complete graph of order $2n+1$;
  \item $G=K_{2n}-I$ is the complete graph of order $2n$ with the edges of a $1$-factor $I$ removed;
  \item $G=K_{2n}+J$ is the multigraph obtained by adding to $K_{2n}$ 
        the edges of a $1$-factor $J$;
  \item $G=\lambda K_v$ is the $\lambda$-fold complete graph.
\end{enumerate}
The most studied problem is represented by $OP(K_v^*; F)$, where 
\[
K_v^*= \left\{
\begin{array}{ll}
K_v & v \mbox{ odd,} \\
K_v - I & v \mbox{ even.}
\end{array}
\right.
\]
In particular, when $v$ is odd, we have the original Oberwolfach problem. The case where $v$ is even
represents the first variant to the original problem, proposed by Huang, Kotzig and Rosa
\cite{HuKoRo79} in 1979. This version of the problem corresponds to finding a maximal packing of $K_v$ ($v$ even) into copies of $F$.
It is known that there is no solution to $OP(K_v^*; F)$
whenever ${F} \in \{[3^2],[3^4],[4,5],[3^2,5]\}$, but that there is a solution for every other instance in which $v \leq 60$ \cite{DFWMR 10, SDTBD}.

When $F$ consists of a single cycle, a solution to $OP(K_v^*; F)$ is a Hamiltonian cycle system of $K_v^*$.  The so-called Walecki construction \cite{Al08} solves this instance of the problem. This construction is notable in that it is 1-rotational ($v$ odd), or 2-pyramidal ($v$-even). That is, it has an automorphism group fixing one or two vertices and acting sharply transitively on the remaining vertices. We discuss such constructions further in Section~\ref{rotational}. 
When $H$ is uniform, $OP(K_v^*; F)$  was solved in the 1990s in
\cite{Alspach Haggvist 85, ASSW, Hoffman Schellenberg 91, RCW}.

In the non-uniform case, OP$(K_v^*; F)$ has been solved if $F$ has exactly two components~\cite{Traetta 13} (see also \cite{BaBrRo97}
for $2$-factorizations of the complete multipartite graph), if ${F}$ is bipartite
(namely, $F$ contains only cycles of even length)~\cite{BryantDanziger} (see Section~\ref{bipirtiteF}), when the order of $F$ belongs to an infinite set of primes \cite{Bryant Schar 09}, or when $F$ has order $2p$ and $p\equiv 5 \pmod{8}$ is a prime  \cite{AlBrHoMaSc 16} (see Section \ref{circulants}).
	Recently, there have been probabilistic results showing eventual existence, but these methods are non constructive and do not provide explicit bounds for existence, so the general problem remains open.

The variant OP($K_{2n}+J; F$) has been formally studied only recently in \cite{BBBSV19},
and it is equivalent to finding a minimum covering of $K_{2n}$ into copies of $F$.
It is completely solved in the uniform case  \cite{AMS87, BBBSV19, LM93}.
In \cite{BBBSV19} the authors point out that the complete solution to OP($K_{2n}-I;F)$ when $F$ is bipartite, mentioned above, implies the solvability of OP$(K_{2n}+J;F)$.

A unifying approach based on constructing $2$-factorizations that have an automorphism group with a prescribed action on the vertex set has recently allowed us to solve completely 
$OP(G;F)$ whenever $F$ contains a sufficiently large cycle and either 
$G=2K_{v}$, or $G= K_{2n+1}, K_{2n}\pm I$ and $F$ has a single-flip automorphism 
(see Section ~\ref{rotational}).
For more results concerning the solvability of infinitely many instances of OP($K_{v}^*; F)$ we refer the reader to the survey \cite[Section VI.12]{Handbook} updated to 2006.

Some recent results concerning OP($\lambda K_v; F$), when $\lambda>1$, can be found in \cite{BaSa21}.
We point out that a solution to $OP(2K_{2n}; F)$ -- which cannot be obtained by doubling a solution of $OP(2K_{2n}, F)$ when $F$ is not bipartite -- together with a solution to $OP(K_{2n}\pm I, F)$ would give rise to a solution of $OP((\lambda K_{2n})\pm I; F)$ for every $\lambda>2$. Therefore, the variant $OP(2K_{2n}, F)$ is equally important as well as all the variants where $\lambda=1$.

Much less is known on the solvability of 
$GOP(K_v^*; \,^{\alpha_1}F_1, \ldots, \,^{\alpha_t}F_t)$ when $t\geq 2$.
Indeed, even the case where $\mathcal{F}$ consists of uniform $2$-factors, that is,
$GOP(K_v^*; \,^{\alpha_1}\ell_1, \ldots, \,^{\alpha_t}\ell_t)$ is still widely open. In \cite{AdaBty06, FraHolRosa, FraRosa00} the problem is solved for odd orders $v$ up to 17, 
and even orders $v$ up to 10 (see also \cite[Sections VI.12.4 and VII.5.4]{Handbook}).
In \cite{BryantDanziger} the problem is settled whenever $v$ is even, each $F_i$ is bipartite (namely, $F_i$ contains only cycles of even length), $\alpha_1\geq 3$ is odd, and the remaining $\alpha_i$ are even. We will consider these types of constructions further in Section~\ref{bipirtiteF}.

The case when the factors have odd cycles seems more difficult. Almost all of the known results are for $t=2$ with uniform factors. The first serious analysis of the case $t=2$ appeared in 2002 \cite{AdaBilBryElz}, which solved existence for uniform cycle sizes $(m, n) = (3, 5), (3, 15)$ and $(5, 15)$. They showed that HWP$(K_{15}, ^6 3, ^1 5)$ does not exist, whilst the existence of HWP$(K_{v}, ^r 3, ^1 5)$ was left open for $v>15$. This highlights the difficulty of the problem when one of the factors only appears once, though this particular case was solved in \cite{WangChenCao}. The case of uniform cycle sizes $(m, n) = (3, 7)$ was almost completely solved in \cite{LeiFu}, with the remaining cases solved in \cite{WangChenCao}. The case $(m, n) = (3, 3x)$ is considered in [5], where the authors solve HWP$(K_v, ^\alpha 3, ^\beta(3x))$ when $v$ is odd and $\beta\neq 1$.

The case of uniform factors with opposite parities seems even harder still. The first result in this case was uniform cycle sizes $(m, n) = (3,4)$, which was almost completely solved in \cite{DanzigerQuattrocchiStevens}. The remaining nine exceptions where given independently in \cite{BonviciniBuratti, OdaOzk17, WangChenCao}. The case of uniform factors with opposite parities was considered more generally in \cite{BDT18b}.
 
These early results have largely been superseded by more general results. In \cite{CaElKhoVan04, ElTipVan02} the problem is solved whenever $v=p^n$ with $p$ a prime number, $t=n$,  and $F_i=[\,^{p^{n-i}} p^i]$, except possibly when $p$ is odd and the first non-zero integer of $(\alpha_1, \alpha_2, \ldots, \alpha_n)$ is $1$. 
Another recent result \cite{BDT19} solves GOP($K_v; \,^{\alpha_1}\ell_1, \ldots, \,^{\alpha_t}\ell_t$)
whenever $v> N=lcm (\ell_1, \ldots, \ell_t)$, $\gcd(\ell_1, \ldots, \ell_t)>1$, 
each $\alpha_i\neq 1$ and they are not all equal to $(N-1)/2$. 
Further results covering specific cases can be found in \cite{Bry04, OdaOzk17, Ste03}.

Many more existence results are known when $t=2$, that is, for HWP$(K_v^*;$ $\,^{\alpha_1}F_1,
\,^{\alpha_2}F_2)$ although they mostly concern the uniform case, which is nonetheless still open. 
Some recent methods that allowed progress in the uniform case \cite{BDT17, BDT18, BDT18b, BDT19b, BDPT23+} will be discussed in Section \ref{blow-up}. 

In the following, we focus on some recent constructive methods that allowed relevant progress to be made, mainly on OP$(G; F)$ whenever $G=2K_v,$ $K_{2n+1}$ or $K_{2n}\pm I$,
on HWP$(K_v^*; \,^{\alpha_1}\ell_1, \,^{\alpha_2}\ell_2)$ and
GOP($K_v; \,^{\alpha_1}\ell_1, \ldots, \,^{\alpha_t}\ell_t$).

\section{Constructive Methods}
\subsection{The uniform case and blow-up type constructions} 
\label{blow-up}

The original problem $OP(K_{2n+1};F)$ and its even variants, $OP(K_{2n}-I; F)$ and $OP(K_{2n}+I;F)$,
have been completely solved in the uniform case, that is, when $F$ contains cycles of the same length. Surprisingly, the uniform case of the Generalized Oberwolfach problem remains widely open, even in the most studied case, that is, the Hamilton-Waterloo problem. 
Nonetheless, the case where all factors of a uniform $2$-factorization are bipartite is peculiar in view of
two general results, Theorems~\ref{BD} and~\ref{refinement}, which are described in Section \ref{bipirtiteF}. Indeed, those results in the uniform case yield the following.

\begin{theorem}[\cite{BryantDanziger, BryantDanzigerDean, Haggkvist 85}] \label{HWP_uniform_even_even}
There is a solution to $\hwp(K_v-I; \,^{\alpha} 2m, \,^{\beta} 2n)$ if and only if $2m$ and $2n$ are both divisors of $v$ and $\alpha+\beta=\frac{v-2}{2}$, except possibly when
\begin{enumerate}
\item $v \equiv 0$ (mod 4), $m \nmid n$, and $1 \in \{\alpha, \beta\}$;
\item $v \equiv 2$ (mod 4), $m \nmid n$, and $\alpha$ and $\beta$ are both odd.
\end{enumerate}
\end{theorem}

In this section, we highlight a few methods that have been recently effectively used to make progress on the uniform GOP.

Let $\Gamma$ be a group, and let $S \subset \Gamma$.
We denote by
$C_g[\Gamma, S]$ ($g\geq 3$) the graph with point set 
$\mathbb{Z}_g\times \Gamma$ and edges $(i,x)(i+1,d+x)$, $i\in \mathbb{Z}_g$, $x\in\Gamma$ and $d\in S$. In other words, 
$C_g[\Gamma,S]$ is the Cayley graph ${\rm cay}(\Z_g\times \Gamma, \{1\}\times S)$ over
 $\Z_g\times \Gamma$ with connection set $\{1\}\times S$; hence, it is $2|S|$-regular.
It is straightforward to see that if $\Gamma$ has order $n$, then 
$C_{g}[\Gamma,\Gamma]$ is isomorphic to the graph $C_g[n]$ obtained as the wreath product of a $g$-cycle $C_g$ by the empty graph of order $n$. In other words, 
$C_g[n]$ is a ``blown-up cycle'' obtained from $C_g$ by replacing every vertex in the cycle with $n$ copies of it: two vertices are adjacent if and only if they come from adjacent vertices of $C_g$. 
Clearly, $C_{g}[\Gamma, S]$ is a subgraph of $C_g[n]$.

A $C_g[n]$-factorization of the complete equipartite graph $K_m[N]$ with $m$ parts, each of size $N$, can be easily constructed by exploiting the following result that settles the existence problem of a $C_g$-factorization of $K_m[z]$.

\begin{theorem}[\cite{Liu00, Liu03}]
\label{Liu} 
Let $g, m$ and $z$ be positive integers with $g\geq 3$.
There exists a $C_g$-factorization of $K_m[z]$ if and only if $g\mid mz$, $(m-1)z$ is even, 
further $\ell$ is even when $m = 2$, and
$(g, m, z) \not\in \{(3, 3, 2), (3, 6, 2), (3, 3, 6),$ $(6, 2, 6)\}$.
\end{theorem}
By blowing up a $C_g$-factorization, we obtain the following straightforward corollary.
\begin{cor}
\label{LiuGen}
Given four positive integers $g, m, n,$ and $z$ with $g\geq 3$,
there exists a $C_g[n]$-factorization of $K_m[nz]$ whenever 
$g\mid mz$, $(m-1)z$ is even, $g$ is even when $t = 2$, and
$(g, m, z) \not\in \{(3, 3, 2)$, $(3, 6, 2), (3, 3, 6), (6, 2, 6)\}$.
\end{cor}
Recall that OP($K_v^*; [^{v/\ell}\ell]$) is completely solved, hence there exists a uniform $2$-factorizations of $K_{nz}$.
Therefore, in view of Corollary \ref{LiuGen}, any uniform $2$-factorization of $C_g[n]$ 
yields a uniform $2$-factorization of $K_{mnz}$, that is, a solution to an instance
of GOP$(K_{mnz}; \,^{\alpha_1}\ell_1, \ldots, \,^{\alpha_s}\ell_s)$.

The approach just described, is then based on constructing solutions to 
$\mathrm{GOP}(C_{g}[\Gamma, S]; \,^{\alpha_1}\ell_1, \ldots, \,^{\alpha_s}\ell_s))$.  Indeed, this approach was integral to the resolution of the uniform Oberwolfach Problem 
OP($K_v^*; [^{v/\ell}\ell]$). Recently, it has also been used to make significant progress towards a complete solution of the uniform Hamilton-Waterloo Problem.

\subsubsection{Projections}

To decompose $C_g[n]$ into $n$-cycles, Alspach, Stinson, Schellenberg and Wagner~\cite{ASSW} introduced the concept of the {\em projection} of an $n$-cycle onto $C_g[n]$.  In our discussion of projections, we will view $C_g[n]$ as $C_g[\mathbb{Z}_n, \mathbb{Z}_n$]; that is, the vertex set will be $\mathbb{Z}_n \times \mathbb{Z}_n$, and edges have the form $(i,x)(i+1,d+x)$, where $i \in \mathbb{Z}_g$ and $x,d \in \mathbb{Z}_n$.  For brevity, we will denote the vertex $(i,x)$ by $i_x$.  We say that the edge $i_{x} (i+1)_{d+x}$ of $C_g[n]$ has {\em difference} $d$.

\begin{definition}
Let $g$ and $n$ be odd integers with $3 \leq g \leq n$, and let $H=(h_0, h_1, \ldots, h_{n-1})$ be a Hamilton cycle of the complete symmetric digraph $\overrightarrow{K}_n$.  The {\em $0$-projection} of $H$ onto $C_g[n]$ is the $n$-cycle
\[
P_0=(0_{h_0}, 1_{h_1}, \ldots, {g-1}_{h_{g-1}}, 0_{h_m}, 1_{h_{m+1}}, 0_{h_{m+2}}, 1_{h_{m+3}}, \ldots, 1_{h_{n-1}}).
\]
The {\em $i$-projection} of $H$ onto $C_m[n]$ is the cycle $P_i$ obtained by adding $i$ modulo $g$ to the first coordinate of every vertex in $P_0$ (that is, we develop the cycle in $\mathbb{Z}_g$).  
\end{definition}
Intuitively, the $i$-projection of $H$ onto $C_g[n]$ wraps once around $C_g$ horizontally, and then zig-zags between parts $i$ and $i+1$.  Figure~\ref{projection} shows an example of a projection of a 9-cycle onto $C_5[9]$.  
\begin{center}
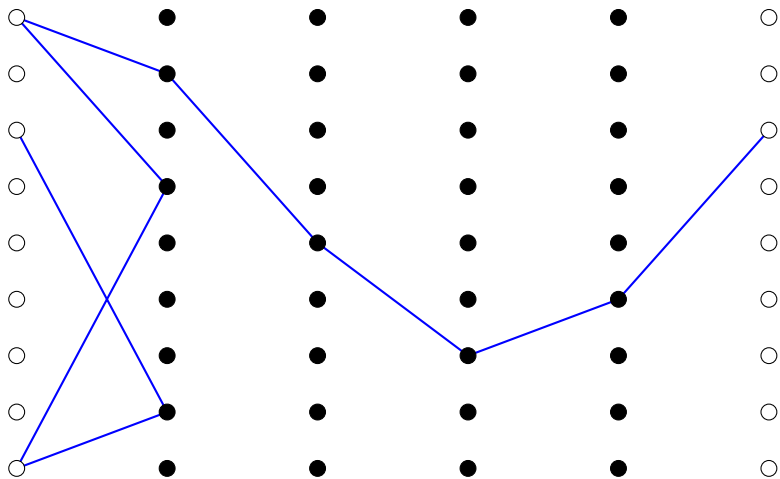
\begin{figure}[h]
\centering
\begin{tikzpicture}[x=2cm,y=0.75cm,scale=1]
\draw[thick, color=blue] (0,8) -- (1,7) -- (2,4) -- (3,2) -- (4,3) -- (5,6);
\draw[thick, color=blue] (0,6) -- (1,1) -- (0,0) -- (1,5) -- (0,8);
\foreach \x in {1,2,3,4}{
	\foreach \y in {0,...,8}{
		\draw[fill=black] (\x,\y) circle (3pt);
	}
}
\foreach \x in {0,5}{
	\foreach \y in {0,...,8}{
		\draw[fill=white] (\x,\y) circle (3pt);
	}
}
\end{tikzpicture}
\caption{The 0-projection of the cycle $(0,1,4,6,5,2,7,8,3)$ on $C_5[9]$.  The columns of white vertices are identified.}
\label{projection}
\end{figure}
\end{center}
The set of $i$-projections of $H$ onto $C_g[n]$ forms a $C_n$-factor $F$ of $C_g[n]$, with the property that for any $x,y \in \mathbb{Z}_n$, $F$ either contains all edges of the form $j_x (j+1)_y$, or else none of these edges.  In particular, if $H$ contains the directed edge $xy$, then $F$ contains either the edges of the form $j_x (j+1)_y$ or $j_y (j+1)_x$, depending on the location of this edge in the cycle.  To account for the other edges, we may use the {\em reverse $i$-projections}, obtained by wrapping around $C_g[m]$ in the opposite order; that is, the reverse $i$-projection of $(h_0, h_1, \ldots, h_{n-1})$ is the cycle
\[
(i_0, (i-1)_{h_1}, \ldots, (i-g+1)_{h_{g-1}}, i_{h_g}, (i-1)_{h_{g+1}}, i_{h_{g+2}}, \ldots, (i-1)_{h_{g+3}}, \ldots, (i-1)_{h_{n-1}}).
\]
While alluded to in~\cite{ASSW}, reverse projections were formalized in~\cite{BDT17}.

To decompose a subgraph of $C_g[n]$ into $C_n$-factors by way of projections, we use a Hamiltonian decomposition of an auxiliary circulant graph Circ$(n;\pm S)$.  Each cycle $H$ in this decomposition is given an arbitrary orientation; the $i$-projections and reverse $i$-projections of $H$ are then taken as 2-factors.  Since the edge $x (x+d)$ of Circ$(n;\pm S)$ will generate all edges of the form $i_x (i+1)_{x_d}$ and $i_{x+d} (i+1)_x$ in these 2-factors, we have the following result.

\begin{theorem}
Let $n$ and $g$ be odd with $3 \leq g \leq n$.  If there is a Hamiltonian decomposition of Circ$(n;\pm S)$, then there is a $C_n$-factorization of $C_g[\mathbb{Z}_n,\pm S]$.
\end{theorem}

Projections alone cannot create 2-factors which contain edges in $C_g[n]$ of difference 0, that is, those of the form $i_x (i+1)_x$.  If needed, $C_n$-factorizations of $C_g[\mathbb{Z}_n, \pm S]$ can be constructed for some set $S$ with $0 \in S$.  In particular, with $S=\{0, 1, 2\}$, we have the following.
\begin{theorem}[\cite{ASSW}]
\label{012}
Let $g$ and $n$ be odd with $3 \leq g \leq n$.  There is a $C_n$-factorization of $C_g[\mathbb{Z}_n, \pm \{0,1,2\}]$.
\end{theorem}
Using this result, combined with the fact that for $n$ an odd prime, every subgraph Circ$(n;\{\pm d\})$ of Circ$(n;\pm \{3,4, \ldots, (n-1)/2\})$ is a Hamiltonian cycle, Alspach, Stinson, Schellenberg and Wagner~\cite{ASSW} conclude that $C_g[n]$ admits a $C_n$-decomposition whenever $g$ is odd, $n$ is an odd prime and $3 \leq g \leq n$.  This result generalizes by also using Hamiltonian decompositions of 4-regular circulants.

\begin{theorem}[\cite{BerFavMah}] \label{DecomposingCirculants}
Every connected 4-regular Cayley graph on a finite Abelian group has a Hamilton cycle decomposition.
\end{theorem}

In particular, every circulant graph of the form Circ$(n;\pm \{d\})$ with $\gcd(d,n)=1$ is trivially decomposable into Hamilton cycles.  Every circulant of the form Circ$(n;\pm \{d_1, d_2\})$ where some linear combination of $d_1$ and $d_2$ is relatively prime with $n$ also 
decomposes into Hamilton cycles by Theorem~\ref{DecomposingCirculants}; in particular this holds when $d_2=d_1+1$.  

These observations lead to the following result. 
\begin{lemma}[\cite{BDT18}]\label{Cn factorizations}
Let $3 \leq g \leq n$ be odd integers and let $S$ be a disjoint union of subsets of $[1,\frac{n-1}{2}]$ of the following forms:
\begin{enumerate}
\item $\{a\}$, where $\gcd(a,n)=1$; 
\item $\{a,b\}$, where $\gcd(xa+yb,n)=1$ for some $x$ and $y$;
\item $[a,b]$, where $b-a+1$ is even;
\item $[a,b] \setminus \{y\}$, where $b-a$ is even and $a \leq y \leq b$;
\item $[a,b] \cup \{x\}$, where $\gcd(x,n)=1$;
\item $([a,b] \setminus \{y\}) \cup \{x\}$, where $\gcd(x,n)=1$;
\item $\{z\} \cup [b,\frac{n-1}{2}]$, where $1 \leq z < b \leq \frac{n-1}{2}$ and $\gcd(b-z,n)=1$.
\end{enumerate}
Then $C_g[\mathbb{Z}_n,\pm S]$ admits a $C_n$-factorization.
\end{lemma}
Note that condition 5 implies that if $S$ is any interval containing an element relatively prime to $n$, then $C_g[\mathbb{Z}_n,\pm S]$ admits a $C_n$-factorization; in particular, taking $S=[3,\frac{n-1}{2}]$ and using Theorem~\ref{012} gives a $C_n$-factorization of $C_g[n]$.

In~\cite{BDT17,BDT18}, the authors use $C_n$-factorizations of $C_m[n]$, together with $C_m$-factorizations derived from row-sum matrices (see Section~\ref{Row-Sum}), to find solutions of the Hamilton-Waterloo Problem when $m$ and $n$ are odd.  In particular,~\cite{BDT17} solves HWP$(K_v;{}^{\alpha}m, {}^{\beta}n)$ when $m$ and $n$ are relatively prime, except possibly if $v=mn$ or $\beta \in \{1,3\}$ and~\cite{BDT18} solves HWP$(K_v;{}^{\alpha}m, {}^{\beta}n)$ for odd $m$ and $n$, except possibly if $\alpha=1$, $\beta \in \{1,3\}$ or $v=\mbox{lcm}(m,n)$ (c.f.\ Theorem~\ref{HWP_uniform}).  Similar results are given in~\cite{BDT18} for $v$ even, but with a few additional possible exceptions.

\subsubsection{Row-Sum Matrices}
\label{Row-Sum}

For this purpose, the concept of a \emph{row-sum matrix} turns out to be very helpful, as described in the following (Theorem \ref{row-sum matrix to factorizations}). We point out that although this concept has been
formally introduced in \cite{BDPT23+}, this type of matrices has been actually used much earlier to build uniform $2$-factorizations $C_{g}[\Gamma, S]$ (see, for example, \cite{BDT17, BDT18,BDT18b,BDT19,BDT19b}).

\begin{definition}
Let $\Gamma$ be a group, and let $S \subset \Gamma$. Also, let $\Sigma$ be an $|S|$-list of elements of $\Gamma$.
A \emph{row-sum matrix} $RSM_\Gamma(S, g; \Sigma)$ is an
$|S|\times g$ matrix, whose $g\geq 2$ columns are permutations of $S$ and such that the list of (left-to-right) row-sums is $\Sigma$. We write $RSM_\Gamma(S, g; \omega(\Sigma))$ whenever we are just interested in the list $\omega(\Sigma)$ of orders of the row-sums.
\end{definition}
The following result, proven in \cite[Theorem 2.1]{{BDT19}} when $\Gamma$ has odd order, 
shows that row-sum matrices can be used to build uniform $2$-factorizations of $C_g[\Gamma, S]$. In \cite{BDPT23+}, the authors point out that the proof of the even case is identical to the odd one.
\begin{theorem}\label{row-sum matrix to factorizations}
If there exists an $RSM_\Gamma(S,g;\Sigma)$, then
$\mathrm{GOP}(C_{g}[\Gamma, S]; g\omega(\Sigma))$ has a solution.
\end{theorem}

The following theorem summarizes the results in \cite{BDT18, BDT18b, BDT19b, KPdiffpar},
which deal with instances of HWP not covered by Theorem \ref{HWP_uniform_even_even}.
\begin{theorem}[\cite{BDT18, BDT18b,BDT19b,KPdiffpar}]\label{HWP_uniform}
There is a solution to HWP$(K^*_v;\,^\alpha M, \,^\beta N)$ when 
\begin{enumerate}
\item $3\leq M < N$ are odd, and $MN/\gcd(M,N)$ divides  $v$,
except possibly if $v= MNu/\gcd(M,N)$ for $u\in\{1,2,4,6\}$, or $\alpha=1$, or $\beta\in\{1,3\}$, or $v$ is even and $(M,N,\beta)=(5,7,5)$;
\item $4\leq M<N$ are even, $M$ does not divide $N$, $\alpha$ and $\beta$ are odd,
except possibly if $1 \in \{\alpha,\beta\}$;
if $\beta=3$, $v\equiv 2 \pmod{4}$ and $\gcd(M,N)=2$;
or if $v=MN/\gcd(M,N) \equiv 2\pmod{4}$, and $\alpha$ and $\beta$ are odd;
\item $N=2^kn$ with $k\geq 1$, $M$ and $n$ are odd, and either $M$ divides $n$, $v>6N>36M$ and $s\geq 3$; or 
$\gcd(M,n)\geq 3$, $4^k$ divides $v$, $v/(4^k\,\mathrm{lcm}(M,n))$ is at least $3$ and $1\not\in\{\alpha,\beta\}$.
\end{enumerate}
\end{theorem}

The following result makes progress on the uniform GOP of odd order.

\begin{theorem}[\cite{BDT19b}]\label{GOP_odd_uniform}
Let $v\geq 3$ be odd, let $3\leq N_1<\ldots< N_t$ and set 
$N=\mathrm{lcm}(N_1,\ldots, N_t)$ and
$g=\mathrm{gcd}(N_1,\ldots, N_t)$; also, let $\alpha_1, \ldots,\alpha_t$ be positive integers.
Then, GOP$(K_v; \,^{\alpha_1}N_1,\ldots, \,^{\alpha_t}N_t)$ has a solution if and only if $N$ is a divisor of $v$ and 
$\sum_{i=1}^t \alpha_i = \frac{v-1}{2}$ except possibly when $t>1$ and at least one of the following conditions is satisfied:
\begin{enumerate} 
\item $\alpha_i= 1$ for some $i\in\{1,\ldots,t\}$;
\item $\alpha_i\in [2, \frac{N-3}{2}]\cup\{\frac{N+1}{2}\}$ for every 
$i\in\{1,\ldots,t\}$;
\item $g= 1$; 
\item $v=N$.
\end{enumerate}
\end{theorem}

Theorems \ref{HWP_uniform} and \ref{GOP_odd_uniform} are essentially obtained by constructing
suitable row sum matrices over abelian groups. The following recent result, which provide further progress on the uniform HWP, is instead obtained also by working over generalized dihedral groups. 

\begin{theorem}\cite{BDPT23+}
Let $v$, $M$ and $N$ be integers greater than 3, and let $\ell=\mbox{lcm}(M,N)$.
A solution to $\mathrm{HWP}(K^*_v;\,^\alpha M, \,^\beta N)$ exists if and only if $M\mid v$ and $N\mid v$, except possibly when
\begin{enumerate} 
\item $\gcd(M,N)\in \{1,2\}$;
\item $4$ does not divide $v/\ell$;
\item $v/4\ell \in\{1,2\}$;
\item $v=16\ell$ and $\gcd(M,N)$ is odd;
\item $v=24\ell$ and $\gcd(M,N)=3$.
\end{enumerate}
\end{theorem}

\subsection{Solutions for bipartite $2$-factors}
\label{bipirtiteF}

In this section we consider the Generalised Oberwolfach problem in the case when the factors in $\cF$ are all bipartite and hence consist of even cycles. Much of this work is based on an observation of H\"aggkvist, which is summarized in the following Lemma.
\begin{lemma}[H\"aggkvist's Doubling Construction \cite{Haggkvist 85}]
	\label{Hagg Lemma}
	Given any $n$-cycle $C_n$, and any bipartite 2-factor of order $2n$, $F$ say, $C_n[2]$ can be factored into two copies of $F$. 
\end{lemma}
The method is illustrated in the picture below.

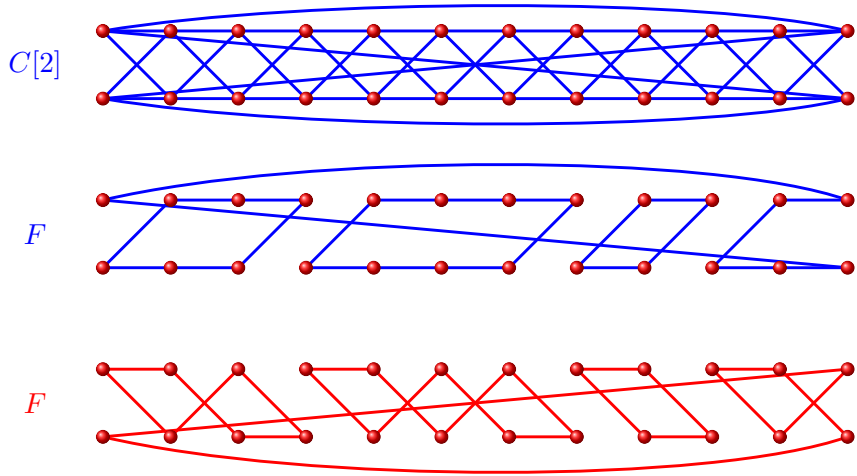
\begin{center}\begin{tikzpicture}[x=1cm,y=1cm,scale=0.9]
	\draw [line width=1.1pt, color=blue] (0,3) -- (11,3);
	\draw [line width=1.1pt, color=blue] (0,4) -- (11,4);
	\draw [line width=1.1pt, color=blue] (0,3) -- (11,4);
	\draw [line width=1.1pt, color=blue] (0,4) -- (11,3);
	\draw [line width=1.1pt, color=blue] (0,3) .. controls (3,2.5) and (9,2.5) ..(11,3);
	\draw [line width=1.1pt, color=blue] (0,4) .. controls (3,4.5) and (9,4.5) .. (11,4);
	\foreach \x in {0,1,...,10} {
		\draw [line width=1.1pt, color=blue] (\x,3) -- (\x+1,4);
		\draw [line width=1.1pt, color=blue] (\x,4) -- (\x+1,3);
	}
	\foreach \x in {0,1,...,11} \foreach \y in {3,4}{
		\shade[ball color=red] (\x,\y) circle (3pt);
	}
	
	\draw [line width=1.1pt, color=blue] (0,.5) -- (2,.5);
	\draw [line width=1.1pt, color=blue] (1,1.5) -- (3,1.5);
	\foreach \x in {0,2}{
		\draw [line width=1.1pt, color=blue] (\x,.5) -- (\x+1, 1.5);
	}
	\draw [line width=1.1pt, color=blue] (3,.5) -- (6,.5);
	\draw [line width=1.1pt, color=blue] (4,1.5) -- (7,1.5);
	\foreach \x in {3,6}{
		\draw [line width=1.1pt, color=blue] (\x,.5) -- (\x+1, 1.5);
	}
	\draw [line width=1.1pt, color=blue] (7,.5) -- (8,.5);
	\draw [line width=1.1pt, color=blue] (8,1.5) -- (9,1.5);
	\foreach \x in {7,8}{
		\draw [line width=1.1pt, color=blue] (\x,.5) -- (\x+1, 1.5);
	}
	\draw [line width=1.1pt, color=blue] (0,1.5) .. controls (3,2.2) and (9,2.2) .. (11,1.5);
	\draw [line width=1.1pt, color=blue] (0,1.5) -- (11,.5);
	\draw [line width=1.1pt, color=blue] (9,.5) -- (11,.5);
	\draw [line width=1.1pt, color=blue] (10,1.5) -- (11, 1.5);
	\draw [line width=1.1pt, color=blue] (9,.5) -- (10, 1.5);

	\draw [line width=1.1pt, color=red] (0,-1) -- (1,-1);
	\draw [line width=1.1pt, color=red] (0,-1) -- (1,-2);
	\foreach \x in {3,7,9}{
		\draw [line width=1.1pt, color=red] (\x,-1) -- (\x+1,-1);
		\draw [line width=1.1pt, color=red] (\x-1,-2) -- (\x,-2);
	}
	\foreach \x in {1,4,5,10}{
		\draw [line width=1.1pt, color=red] (\x,-1) -- (\x+1, -2);
		\draw [line width=1.1pt, color=red] (\x,-2) -- (\x+1, -1);
	}
	\foreach \x in {2,3,6,7,8,9}{
		\draw [line width=1.1pt, color=red] (\x,-1) -- (\x+1, -2);
	}
	\draw [line width=1.1pt, color=red] (0,-2) .. controls (3,-2.7) and (9,-2.7) .. (11,-2);
	\draw [line width=1.1pt, color=red] (0,-2) -- (11,-1);
	\
	\draw (-1,3.5) node[blue]{$C[2]$};
	\draw (-1,1) node[blue] {$F$};
	\draw (-1,-1.5) node[red] {$F$};
	\foreach \x in {0,1,...,11} \foreach \y in {.5,1.5}{
		\shade[ball color=red] (\x,\y) circle (3pt);
	}
	\foreach \x in {0,1,...,11} \foreach \y in {-1,-2}{
		\shade[ball color=red] (\x,\y) circle (3pt);
	}
	\end{tikzpicture}
\end{center}

As noted earlier, the classical Walecki construction, see for example \cite{Al08}, provides a factorization of $K_v$ into Hamiltonian cycles, $H_i$ say, when $v$ is odd. Each of these can then be ``doubled" to $H_i[2]$, these can then be factored into pairs of isomorphic bipartite factors by the previous lemma. This yields the following result.

\begin{lemma}[\cite{Haggkvist 85}]
Let $v\equiv 2\bmod 4$ and let $\mathcal{F}=\{^{\alpha_1}F_1, ^{\alpha_2}F_2, \ldots, ^{\alpha_t}F_t\}$ be a collection of bipartite $2$-regular graphs of order $v$, where $\alpha_1, \alpha_2, \ldots, \alpha_t$ are nonnegative even integers satisfying $\alpha_1+\alpha_2+\ldots+\alpha_t = \frac{v-1}{2}$, then GOP$(K_v; \mathcal{F})$ has a solution.
\end{lemma}

More recently, Bryant and Danziger \cite{BryantDanziger} have adapted this method to the case $v \equiv 0\bmod 4$. In this case they show that $K_v$ can be factored into a specific 3-regular circulant graph, $G$, and $\frac{v-4}{2} $ Hamiltonian factors.
More specifically, when $v$ is even, given a difference $d\in \Z_v$, we denote the set of edges with difference $d$ starting from an even vertex by $d^{\,\rm e}$ and the set of edges with difference $d$ starting from an even vertex by $d^{\,\rm o}$.
In \cite{BryantDanziger} it is shown that for every $m\geq 8$ there is a factorization of $K_v$ into $\frac{v-4}{2} $ Hamiltonian factors and a copy of the 3-regular circulant containing all edges of the form $\{1,3^{\,\rm e}\}$, denoted $\langle1,3^{\,\rm e}\rangle$. This decomposition is then doubled and the Hamilton factors factored using Lemma~\ref{Hagg Lemma}. They then show that the remaining 7-regular graph (including each vertex $x$ and its doubled image as an edge) can be factored into three copies of any bipartite factor. This results in the following general result.

\begin{theorem}[\cite{BryantDanziger, Haggkvist 85}]\label{BD}
	 Let $\mathcal{F}=\{^{\alpha_1}F_1, ^{\alpha_2}F_2, \ldots, ^{\alpha_t}F_t\}$ be a collection of bipartite $2$-regular graphs of order $v$, where  $\alpha_1+\alpha_2+\ldots+\alpha_t = \frac{v-2}{2}$ and  $\alpha_2, \ldots, \alpha_t$ are nonnegative even integers. 
	 If $\alpha_1$ is nonnegative and \\
	 \begin{tabular}{ll}
	 	when $v \equiv 2\bmod 4$, $\alpha_1$ is even, or  \\
	 	when  $v \equiv 0\bmod 4$, $\alpha_1\geq 3$ is odd,
	 \end{tabular}\\
	then GOP$(K_v; \mathcal{F})$ has a solution.
\end{theorem}
This theorem then solves the Oberwolfach problem OP$(K_v-I,F)$, when $v$ is even and $F$ is a bipartite factor. 

In \cite{BryantDanzigerDean}, Bryant, Danziger and Dean extended this method to consider the Hamilton-Waterloo problem in particular. They showed that the doubled 7-regular graph $\langle 1,3^{\,\rm e}\rangle[2]$ could be factored into three copies of a bipartite factor, $F$, and a 1-factor consisting of the edges $x$ and its double, except possibly when $F \in \{[6^r], [4, 6^r ] \mid r \equiv 2 \bmod 4\}$. 

Given two bipartite factors $F_1$ and $F_2$ we say that $F_1$ is a {\em refinement} of $F_2$, if $F_1$ can be obtained from $F_2$ by replacing each cycle of $F_2$ with a bipartite 2-regular graph on the same vertex set. For
example (see \cite{BryantDanzigerDean}), $[4, 8^3, 10^2 , 12]$ is a refinement of $[4, 16, 18, 22]$, but $[4, 18^2, 20]$ is not. They where able to give a complete solution of the Hamilton-Waterloo problem whenever $F_1$ is a refinement of $F_2$.
\begin{theorem}
	[\cite{BryantDanzigerDean}] \label{refinement}
	Let $F_1$ and $F_2$ be bipartite $2$-regular graphs of order $v$ such that $F_1$ is a refinement of $F_2$.  
	There is a solution to $\hwp(K_v-I; \,^{\alpha}F_1, \,^{\beta} F_2)$ 
	if and only if $\alpha+\beta=\frac{v-2}{2}$.
\end{theorem}

Bryant Danziger and Pettersson~\cite{BryantDanzigerPettersson} applied these methods to multipartite graphs $K_r[n]$, and where able to completely solve OP$(K_r[n], F)$ for bipartite factors $F$.
\begin{theorem}
	If $r\geq 2$ and $F$ is a bipartite 2-regular graph of order $rn$, then OP$(K_r[n], F)$ has a solution if and only if $n$ is even; except that there is no 2-factorization of $K_2[6]$ $(\cong K_{6,6})$ into $[6, 6]$.
\end{theorem}

\subsection{2-Factorizations of circulant graphs}\label{circulants}

A major breakthrough in the Oberwolfach Problem came in a 2009 paper by Bryant and Scharaschkin~\cite{Bryant Schar 09}, which gave the first complete solution for every 2-factor on an infinite family of orders $v$.  In particular, the following theorems are the main results of~\cite{Bryant Schar 09}.
\begin{theorem}[\cite{Bryant Schar 09}] \label{BS1}
For $v \equiv 1$ (mod 16), let $G_v$ and $H_v$ be the multiplicative subgroups of $\mathbb{Z}_v$ of index 2 and 8, respectively.  If $v \equiv 1$ (mod 16) is prime, $\{1,2,3,4\} \subseteq G_v$ and $G_v/H_v=\{H_v,2H_v,3H_v,4H_v\}$, then there is a solution of OP$(K_v;F)$ for every $2$-factor $F$ of order $v$.  Moreover, there are infinitely many such $v$.
\end{theorem}

\begin{theorem}[\cite{Bryant Schar 09}] \label{BS2}
Let $p \equiv 5$ (mod $24$) be prime.  There is a solution of OP$(K_v-I; F)$ for each 
$2$-factor $F$ of order $v=2p$.  
\end{theorem}

The fact that Theorem~\ref{BS2} gives a complete solution for an infinite family of orders is a direct consequence of Dirichlet's Theorem on primes in arithmetic progression; the infinitude of the orders covered by Theorem~\ref{BS1} comes from an application of the Chebotarev Density Theorem from algebraic number theory.

The authors of~\cite{AlBrHoMaSc 16} gave an improvement of Theorem~\ref{BS2}.
\begin{theorem}[\cite{AlBrHoMaSc 16}] \label{ABHMS}
Let $p \equiv 5$ (mod $8$) be prime.  There is a solution of OP$(K_v-I; F)$ for each 
$2$-factor $F$ of order $v=2p$.  
\end{theorem}

The solutions to the Oberwolfach Problem constructed by Theorems~\ref{BS1}, \ref{BS2} and~\ref{ABHMS} follow the same strategy:  first factor $K_v$ or $K_v-I$ into circulant graphs of the form Circ$(n;\pm\{1,2,3,4\})$, and then construct appropriate factorizations of this circulant.  In~\cite{AlBrHoMaSc 16}, some additional results are obtained on the Oberwolfach Problem using 2-factorizations of other circulant graphs, namely factorizations of Circ$(n;\pm\{1,2\})$ into a Hamiltonian cycle and a copy of a given 2-regular graph $F$ from~\cite{Bry04, Rodger 04}, factorizations of Circ$(n;\pm\{1,2,3\})$ into 2-factors with a limited number of 3-cycles, and factorizations of Circ$(n;\pm\{1,3,4\})$ into 2-factors whose minimum cycle length is at least 6.

\subsection{Solutions with well-behaved symmetries}
\label{rotational}
Constructing $2$-factorizations with well-behaved automorphisms has proved to be a winning method for solving a large portion of the original Oberwolfach problem and some of its variants. 
A $2$-factorization $\mathcal{F}$ of a graph $G$ is called 
\begin{enumerate}
  \item \emph{$r$-rotational} ($r\geq1$) over $\Gamma$ 
  if $\Gamma$ is an automorphism group of $\mathcal{F}$ that fixes one vertex and 
  generates on the remaining vertices $r$ $\Gamma$-orbits, each of size $|\Gamma|$;
  \item \emph{$f$-pyramidal} ($f\geq0$) over  $\Gamma$ if $\Gamma$ is an automorphism group of
  $\mathcal{F}$ fixing $f$ vertices and acting sharply transitively on the remaining.
\end{enumerate}
Note that $1$-rotational is equivalent to $1$-pyramidal, while $0$-pyramidal solutions are usually referred to as \emph{sharply-transitive} solutions. Note also that given a group $\Gamma$ with a subgroup $\Sigma$, any $1$-rotational solution to OP over $\Gamma$ is also $r$-rotational over $\Sigma$, with $r=|\Gamma|/|\Sigma|$. The approach based on building $r$-rotational $2$-factorizations has been successfully used to solve all instances of 
OP$(K^*_{v};F)$ for all orders $v\leq40$ \cite{DFWMR 10} and whenever $41\leq v\leq 60$ \cite{SDTBD}.

The structure of a sharply transitive $2$-factorization of $K_v$ was studied in 
\cite{Buratti Rinaldi 05}. 
Considering that regular $2$-factorizations of $K_v$ have been built mainly over the cyclic group and in the uniform case (see, for example, \cite{BuRaZu05,  Buratti Rinaldi 05}), this suggests that sharply-transitive solutions to OP are quite rare. Structural results on $f$-pyramidal $2$-factorizations of $K_v$ can be found in \cite{BMR09}. Similarly, the structure of $2$-pyramidal $2$-factorizations of $K_{2n}-I$ has been studied in \cite{BuTr15}.

In the following, we focus on the $1$-rotational method and the reason is twofold.
 First of all, it is known that every 1-rotational solution to OP$(F)$ of order $2n+1$ gives rise to a $2$-pyramidal \cite{BuTr15, HuKoRo79} (resp.\  sharply-transitive \cite{BuDaTr22}) solution to OP$(K_{2n+2}-I; F^*)$ (resp.\  OP$(K_{2n}+I;F^*)$)
for a suitable $2$-factor $F^*$ of order $2n+2$ (resp. $2n$).
Also, the $1$-rotational method has lately allowed a large portion of the original problem to be solved,
and (based on the previous observation) its maximum packing and minimum covering variants.

Solutions to OP that are $1$-rotational over an arbitrary group $G$ were formally investigated in 2008 \cite{Buratti Rinaldi 08}, although this type of solution has appeared earlier, mainly over the cyclic group \cite{BuZu98, BuZu01}, and in some cases making use of \emph{terraces} in groups \cite{OlPr03, Ol05}. 
It is worth pointing out that the well-known Walecki construction yields a $1$-rotational solution to OP$(K_{2n+1};[2n+1])$ (see, for example, \cite{Al08}). 

The existence of a $1$-rotational solution to $OP(K_{2n+1};H)$ and $OP(2K_{v}; H)$ over a group $G$ is shown in 
\cite{Buratti Rinaldi 08, Buratti Traetta 12} to be equivalent to the existence of a copy of $H$ with vertices
in $G\,\cup\,\{\infty\}$ satisfying further algebraic conditions, which we describe in the following after introducing the required terminology and notation.

Let $G$  be an arbitrary group and let $\Gamma$ be a graph with vertices in 
$G \,\cup\, I$, where  $I= \{\infty_1, \ldots, \infty_f\}$ is a set disjoint from $G$.
Also,  for every $g \in G$, we denote by $\Gamma+g$ \emph{the translate of $\Gamma$ by $g$}, that is,  the graph obtained from $\Gamma$ by replacing each vertex $x\in V(\Gamma)\setminus I$ with $x+g$. We also denote by $\Delta \Gamma = \{x-y\mid (x,y)\in G^2, \{x,y\}\in E(\Gamma)\}$ \emph{the list of differences of $\Gamma$}, that is, the list of all differences between adjacent vertices of $\Gamma$ not belonging to $I$. \emph{The $G$-development of $\Gamma$} is the multiset $Dev_G(\Gamma)$ of all translates of $\Gamma$. Similarly, \emph{the $G$-orbit of $\Gamma$} is the multiset $Orb_G(\Gamma)$ of all `distinct' translates of $\Gamma$ by all the elements of $G$.

\begin{defn}[\cite{Buratti Rinaldi 08, Buratti Traetta 12}] Let $F$ be a $2$-regular graph with vertex set $V(F)=G\,\cup\, \{\infty\}$. We say that $F$ is a \emph{twofold $2$-starter of $G$} if 
\begin{enumerate}
  \item $\Delta F = \,^2{(G\setminus\{0\})}$, that is, $\Delta F$ contains exactly every non-zero element of $G$ with multiplicity $2$.
\end{enumerate}
Also, if $F$ satisfies the following further condition
\begin{enumerate}
  \setcounter{enumi}{1}
  \item there exists an element  $y\in G$ of order $2$ (i.e., an involution of $G$) such that $F+y=F$,
\end{enumerate}
we say that $F$ is a \emph{$2$-starter of $G$}.
\end{defn}

These concepts were introduced in \cite{Buratti Rinaldi 08, Buratti Traetta 12}  to characterize $1$-rotational solutions to OP. More precisely, we have the following.

\begin{theorem}[\cite{Buratti Rinaldi 08, Buratti Traetta 12}]\label{rotationalstarter}
  \begin{enumerate}
  \item 
  $OP(K_{2n+1};H)$ has a $1$-rotational solution over a group $G$ if and only if 
  $G$ has a $2$-starter $F$ isomorphic to $H$. In this case, $Orb_G(F)$ is 
  a $1$-rotational solution to OP$(H)$. 
  \item 
  $OP(2K_v; H)$ has a $1$-rotational solution over a group $G$ if and only if 
  $G$ has a twofold $2$-starter $F$ isomorphic to $H$. 
  In this case, $Dev_\Gamma(F)$ is 
  a $1$-rotational solution to OP$(2K_v; H)$. 
  \end{enumerate}
\end{theorem}
The above characterization translates the existence problem of $1$-rotational $2$-factorizations into the construction of suitable starters over an arbitrary group.
This approach has proven successful (since 
it was formally introduced  in \cite{Buratti Rinaldi 08}) 
to construct solutions to OP$(K_v;F)$.
In \cite{BuRi09} the authors completely characterize the $1$-rotational solutions of 
OP$(K_{2s+3};[3,2s])$ and via a composition technique infinitely many $1$--rotational solutions are built in \cite{RiTr11}.
Other $1$-rotational constructions include those given in \cite{BuZu98, BuZu01} concerning $OP(K_{6n+3};[\,^{2n+1}3])$, 
and those given in a series of papers
\cite{OlPr03, Ol05, OlSt09, OlWil11} which mainly focus on the case where $F$ has three components.

In \cite{Buratti Traetta 12}, the authors notice that \emph{graceful labellings} can be used to  build twofold $2$-starters over the cyclic group. More precisely, assume there exists a graceful labelling of a $2$-regular graph $H$ of order $v$, that is,  a copy $F$ of $H$ such that
\[
  V(F)=\{0,\ldots, v-1\}\,\cup\, \{\infty\}\;\;\;\text{and}\;\;\; 
  \Delta F=\pm \{1,\ldots, v-1\}.
\]
Clearly, reviewing the vertices of $F$ modulo $v$ provides a twofold $2$-starter of
$\Z_v$, hence (by Theorem \ref{rotationalstarter}) a solution to $OP(2K_v;H)$.
Graceful labellings were used in \cite{Traetta 13} to settle $OP(2K_v:F)$
with $1$-rotational solutions whenever $F$ has exactly two cycles. 
We point out that when $F$ has two components, a complete solution to $OP(2K_v;F)$ had been given earlier in \cite{BaBrRo97}
using different techniques which do not aim to  build solutions having well-behaved automorphisms.
The same approach based on constructing graceful labellings of suitable $2$-regular graphs was recently used in \cite{BuDaTr22} to completely settle 
OP$(2K_v;F)$ with $1$-rotational solutions whenever $F$ contains a cycle of length greater than an explicit lower bound which depends of the remaining cycle lengths.

\begin{theorem}[\cite{BuDaTr22}]
Let $F=[h, \ell_1,\ldots, \ell_r, \ell'_1,\ldots, \ell'_s]$ be any simple $2$-regular graph,
where the $\ell_i$s are even and the $\ell'_j$s are odd. 
Also, set $\ell=\max\{\ell_1,\ldots, \ell_r\}$ and  $\ell'=\max\{\ell'_1,\ldots, \ell'_s\}$.
Then OP$(2K_v; F)$ has a $1$-rotational solution whenever
\[h> 
  \begin{cases}
    12r(\ell+3) + 7^{s}(2\ell'+1) -6  &\text{if $r,s > 0$}, \\    
    2(r+1)(\ell+1) +5 &\text{if $r> 0=s$},  \\    
    3\cdot7^{s-1}(2\ell'+1) &\text{if $s> 0=r$}. \\       
  \end{cases} 
\]   
\end{theorem} 

A powerful doubling construction, described in \cite{Buratti Traetta 12}, 
builds $2$-starters of $\Z_{2n}$ by doubling 
a twofold $2$-starter of $\Z_{n}$. In view of Theorem \ref{rotationalstarter},
this is equivalent to saying that any $1$-rotational solution to
OP$(2K_{n+1};H)$ over $\Z_n$ yields a $1$-rotational solution to OP$(K_{2n+1};F)$ over $\Z_{2n}$ and for a suitable $2$-regular graph $F$ whose cycle-structure depends on $H$. 

\begin{theorem}[Doubling construction \cite{Buratti Traetta 12}]
\label{doubling}
If there exists a twofold $2$-starter $H\simeq[\ell_\infty,$ $\ell_1,\ldots, \ell_t]$ of $\Z_n$, then there exists a $2$-starter $F$ of $\Z_{2n}$, with
\begin{equation}\label{doubling:structure}
  F\simeq [2\ell_\infty-1] \cup [2\ell_i \mid i \in I] \cup [\,^2\ell_j \mid j \in J],
\end{equation}
for a suitable partition $I\,\cup\, J$ of the set $\{1,2,\ldots,t\}$.
\end{theorem}

As mentioned in \cite{Buratti Traetta 12}, the doubling construction works over any group $G$ with a unique involution $y$: note that the cyclic groups fall into this class.
In fact, the construction can be adapted and applied to any twofold 2-starter of the quotient group $G/\{0,y\}$, yielding this way a suitable $2$-starter of $G$. 
This observation may be helpful to find $1$-rotational solutions to $OP(K_{2n+1};F)$
(or equivalently, $2$-starters) when they cannot be found over the cyclic group.
Without going into details (for which we refer the reader to \cite{Buratti Traetta 12}), further conditions on $H$ in Theorem \ref{doubling} guarantee the existence of
a 2-starter of $\Z_{2n}$ as in \eqref{doubling:structure}, for 
every partition $I\,\cup\, J$ of the set $\{1,2,\ldots,t\}$.

The doubling construction was used in \cite{Traetta 13}  to settle OP$(F)$   with $1$-rotational solutions whenever $F$ has exactly two cycles. Another recent successful application of the doubling construction allowed the authors of the present survey \cite{BuDaTr22} to make relevant progress 
(Theorem \ref{singleflip}) on the ($1$-rotational) solvability of OP$(K_{2n+1};F)$. 
To simplify the statement of our results, we introduce the following notation.
Given a list of positive integers $L=\{\ell_1, \ldots, \ell_s\}$ (not necessarily distinct), we write $L_0$ to represent the multiset of even elements of $L$ greater than 2, and $L_1$ to represent the multiset of odd elements of $L$. Note that by $|L|$ we mean the size of $L$ as
a multiset.

\begin{theorem}\label{singleflip}
Let $F=[h, 2\ell_1, \ldots, 2\ell_r, \;^{2}\ell_{r+1}, \ldots, \;^{2}\ell_{s}]$ 
where
$2\leq \ell_1< \ell_2< \ldots < \ell_r$, 
$\ell_{r+1},\ell_{r+2}, \ldots, \ell_{s} \geq 3$ and $h\geq 3$ is odd.
Also, let
$L=\{\ell_1, \ell_2, \ldots, \ell_s\}$.
Then $OP(K_{2n+1};F)$ has a $1$-rotational solution whenever 
\[h> 16\max(1,h_0) + 20\max(3,h_1) + 29,\] 
where
$h_0 = 2|L_0|\,(\max(L_0) +3) - 1$ and $h_1 = 7^{s-|L_0|-1}(2\max(L_1)+1)$.
\end{theorem}
Without going into further details, since every $1$-rotational solution to the original OP yields solutions to $OP(K_{2n}-I;F)$ and $OP(K_{2n}+I;F)$, a result similar to 
Theorem \ref{singleflip} (see \cite{BuDaTr22}) holds for the maximum packing and the minimum covering version of OP.

\end{document}